\documentclass{amsart}

\usepackage{graphicx}
\usepackage{amsmath,amscd}
\usepackage{amssymb}
\usepackage{amsfonts}
\usepackage{tikz}
\usetikzlibrary{tqft}
\usetikzlibrary{shapes}
\usepackage{cite}
\usepackage{tikz-cd} 

 \newtheorem{theorem}{Theorem}[subsection]
 
 \newtheorem{lemma}[theorem]{Lemma}

 \theoremstyle{definition}
 \newtheorem{definition}[theorem]{Definition}
 \theoremstyle{remark}
  \newtheorem{notation}[]{Notation}
 \numberwithin{equation}{subsection}

\newcommand{\ra}{\rightarrow}

\begin{document}

\title[The Adjoint Reidemeister Torsion for Compact 3-Manifolds Admit a Unique Decomposition]
{The Adjoint Reidemeister Torsion for Compact 3-Manifolds Admit a Unique Decomposition}

\author{Esma Dirican Erdal}

\address{Deptartment of Mathematics,
\.{I}zmir Institute of Technology,
 35430, \.{I}zmir, Turkey}
 
\email{esmadiricanerdal@gmail.com, esmadirican131@gmail.com}

\keywords{Adjoint Reidemeister torsion, Compact $3$-manifolds, Disk sum}


\date{\today}

\dedicatory{}

\commby{}

\begin{abstract}
Let $M$ be a triangulated, oriented, connected compact $3$-manifold with connected non-empty boundary. Such a manifold admits a unique decomposition into
$\triangle$-prime $3$-manifolds. In this paper, we show that the adjoint Reidemeister torsion has a multiplicative property on the disk sum decomposition of compact $3$-manifolds without a corrective term.
\end{abstract}

\maketitle

\section{Introduction}
 
 Along this paper we concern the class of triangulated, oriented, connected compact $3$-manifolds with connected non-empty boundary. Let $\chi$ denote the class of $3$-manifolds $M$ with connected non-empty boundary such that every $2$-sphere in $M$ bounds a $3$-cell. Let $M$ and $M'$ be two manifolds in the class $\chi.$ Then the disk sum (boundary connected sum) $M \bigtriangleup M'$ can be formed by pasting a $2$-cell on the boundary of $M$ to a $2$-cell on the boundary of $M'.$ The operation of disk sum $\bigtriangleup$ is well-defined, associative, and commutative up to homeomorphism. A manifold $M \in \chi$ is called $\triangle$-prime if it is not a $3$-cell, and whenever $M \cong P\triangle P'$, either $P$ or $P'$ is a $3$-cell. In \cite{JLGROSS1}, Gross proved the following decomposition theorem:
\begin{theorem}[\cite{JLGROSS1}]\label{gross1}
 For any $3$-manifold $M$ (different from a $3$-cell) with connected non-empty boundary, there is an isomorphism 
 $$M \cong M_1 \triangle M_2 \triangle \ldots \triangle M_k,$$ 
where the summands $M_i$ are $\triangle$-prime $3$-manifolds, and they are uniquely determined up to order and homeomorphism.
\end{theorem}
 
 In 1935, Reidemeister introduced a new invariant, called now \textit{Reidemeister torsion}, to classify $3$-dimensional lens spaces (up to PL equivalence) \cite{Reidemeister}. Later,
Franz classified higher dimensional lens spaces by extending the notion of this invariant \cite{Franz}. 
In 1969, Kirby and Siebenmann showed that Reidemeister torsion is a topological invariant for manifolds \cite{RCLC}. The invariance for arbitrary simplicial complexes was proved by Chapman \cite{T. A. Chap} and thus the classification of lens spaces of Reidemeister and Franz was shown to be a topological invariant. In 1961, Milnor disproved Hauptvermutung by using this invariant. He constructed two homeomorphic but combinatorially distinct finite simplicial complexes. Then he described Reidemeister torsion with the Alexander polynomial which plays an important role in knot theory and links \cite{Milnor2,Milnor}.

The twisted chain complex of $3$–manifolds $M$ with connected non-empty boundary is never acyclic, so the computation of the adjoint Reidemeister torsion of $M$ involves a corrective term coming from the homologies. In this paper, we show that the adjoint Reidemeister torsion of $M$ has multiplicativity property without a corrective term. More precisely, we establish a multiplicative adjoint Reidemeister torsion formula for $M$ in terms of the adjoint Reidemeister torsions of $\triangle$-prime $3$-manifolds in the unique disk sum decomposition:
  \begin{theorem}\label{theo1}
If $M=\overset{n}{\underset {i=1}{\bigtriangleup}}(M_i)$ and $\psi_{_i}:\pi_1(M_i)\rightarrow G$ is a given representation for each $i,$ then the followings hold:
\begin{itemize}
\item[(i)]{There exists a unique representation $\varrho: \pi_1(M)\rightarrow G,$}
\item[(ii)] {For a given basis $\mathbf{h}^{M}_p$ of $H_p(M;\mathfrak{g}_{\mathrm{Ad}_{\varrho}})$ and a fixed  basis $\mathbf{h}_{0}^{\mathbb{D}^2}=\varphi_0(\mathbf{c}_0)$ of $H_0(\mathbb{D}^2;\mathfrak{g}_{\mathrm{Ad}_{\varrho_{|_{\mathbb{D}^2}}}}),$ there exists a basis 
$\mathbf{h}^{M_i}_p$ of $H_p(M_i;\mathfrak{g}_{\mathrm{Ad}_{\psi_{_{M_i}}}})$  for each $i\in\{1,\ldots,n\}$ such that the following formula is valid
\begin{eqnarray*}
 \mathbb{T}_{\varrho}(M,\{\mathbf{h}^{M}_p\}_{p=0}^3)= \prod_{i=1}^{n}
 \mathbb{T}_{{\psi_{_i}}}(M_i,\{\mathbf{h}_p^{M_i}\}_{p=0}^{3}).
 \end{eqnarray*}}
\end{itemize}
\end{theorem} 

\section{The Adjoint Reidemeister Torsion}

Along this paper $G$ denotes a complex reductive algebraic group $SL_n(\mathbb{C})$ or $PSL_n(\mathbb{C}).$ Let $X$ be a CW-complex with finite dimension $n$ and let $\widetilde{X}$ denote its universal covering. Let us denote the non-degenerate Killing form on $\mathfrak{g}$ by $\mathcal{B}$ which is defined by $\mathcal{B}(A,B) = 4 \;\mathrm{Trace}(AB).$ For a representation $\rho:\pi_1(X)\rightarrow G,$ let us consider the action of $\pi_1(X)$ on $\mathfrak{g}$ via the adjoint of $\rho.$ The integral group ring is defined by
 $$\mathbb{Z}[\pi_1(X)]=\left\{\sum_{i=1}^p m_i\gamma_i\; ; m_i\in \mathbb{Z},\;
\gamma_i\in \pi_1(X),\; p\in \mathbb{N}\right\}.$$

Let $K$ be a cell-decomposition of $X$ and $\widetilde{K}$ be a lifting of $K.$ By using the cellular chain complex $C_{\ast}(\widetilde{K}; \mathbb{Z}),$ one can define the twisted chains as follows 
 \begin{equation}\label{defchn1}
C_{\ast}(K;\mathfrak{g}_{\mathrm{Ad}_\rho}):=\displaystyle
C_{\ast}(\widetilde{K};\mathbb{Z})\displaystyle\otimes \mathfrak{g}
/\sim, \end{equation}
where $\sigma\otimes t \sim
\gamma\cdot\sigma\otimes\gamma\cdot t,\forall \gamma\in \pi_1(X),$
the action of $\pi_1(X)$ on $\widetilde{X} $ is the deck
transformation, and the action of $\pi_1(X)$ on $\mathfrak{g}$ is
the adjoint action.

Let $\{e^p_j\}_{j=1}^{m_p}$ be the generators for $C_p(K;\mathbb{Z}).$
Fixing a lift $\widetilde{e}^p_j$ of $e^p_j$ in $\widetilde{X}$ for each $j\in \{1,\ldots, m_p\},$ we get a
$\mathbb{Z}[\pi_1(X)]$-basis
$c_p=\{\widetilde{e^p_j}\}_{j=1}^{m_p}$ for
$C_p(\widetilde{K};\mathbb{Z}).$ Suppose that
$\mathcal{A}=\{\mathfrak{a}_k\}_{k=1}^{\dim \mathfrak{g}}$ is a
$\mathcal{B}$-orthonormal basis of $\mathfrak{g}$.
Then $\mathbf{c}_p=c_p\otimes_{\rho}
\mathcal{A}$ is a
\emph{geometric basis} for $C_p(K;\mathfrak{g}_{\mathrm{Ad}_
\rho}).$ 

Consider the following chain complex
\begin{equation*}\label{chaincomplex}
  \begin{array}{ccc}
 C_{\ast}:=C_{\ast}(K;\mathfrak{g}_{\mathrm{Ad}_{\rho}})=(0 \to C_{n}(K;\mathfrak{g}_{\mathrm{Ad}_
\rho}) {\rightarrow}
C_{n-1}(K;\mathfrak{g}_{\mathrm{Ad}_
\rho})\rightarrow \cdots \rightarrow C_{0}(K;\mathfrak{g}_{\mathrm{Ad}_
\rho}){\rightarrow}0).
  \end{array}
 \end{equation*}

 For $p\in \{0,\ldots,n\},$ let
$$B_p(C_{\ast})=\mathrm{Im}\{\partial_{p+1}:C_{p+1}\rightarrow C_{p}
\},$$
$$Z_p(C_{\ast})=\mathrm{Ker}\{\partial_{p}:C_{p}\rightarrow
C_{p-1} \},$$ and $H_p(C_{\ast})=Z_p(C_{\ast})/B_p(C_{\ast})$ be
$p$-th homology group of the chain complex $C_{\ast}(K;\mathfrak{g}_{\mathrm{Ad}_{\rho}}).$ Then there are
the following short exact sequences
\begin{equation}\label{Equation1}
0\longrightarrow Z_p(C_\ast) \stackrel{\imath}{\hookrightarrow} C_p(C_{\ast})
\stackrel{\partial_p}{\longrightarrow} B_{p-1}(C_\ast) \longrightarrow 0,
\end{equation}
\begin{equation}\label{Equation2}
0\longrightarrow B_p(C_\ast) \stackrel{\imath}{\hookrightarrow} Z_p(C_\ast)
\stackrel{\varphi_p}{\longrightarrow} H_p(C_\ast) \longrightarrow 0.
\end{equation}
Here, $\imath$ and $\varphi_p$ are the inclusion and the natural
projection, respectively.

Let $s_p:B_{p-1}(C_{\ast})\rightarrow C_p(C_{\ast}),$
$\ell_p:H_p(C_{\ast})\rightarrow Z_p(C_{\ast})$ be sections of
$\partial_p:C_p(C_{\ast}) \rightarrow B_{p-1}(C_{\ast}),$
$\varphi_p:Z_p(C_{\ast})\rightarrow H_p(C_{\ast}),$
respectively. By Splitting Lemma, the short exact sequences (\ref{Equation1})
and (\ref{Equation2}) yield
\begin{equation}\label{Equation0}
C_p(C_{\ast})=B_{p}(C_{\ast})\oplus \ell_p(H_p(C_{\ast }))\oplus s_p(B_{p-1}(C_{\ast})).
\end{equation}
If $\mathbf{b_p}$ and
$\mathbf{h_p}$ are respectively bases of $B_p(C_\ast),$ and $H_p(C_\ast),$ then, by equation~(\ref{Equation0}), the following disjoint union
$$\mathbf{b}_p\sqcup \ell_p(\mathbf{h}_p)\sqcup
s_p(\mathbf{b}_{p-1})$$
becomes a new basis for $C_p.$

\begin{definition}The adjoint Reidemeister torsion of a chain complex $C_{\ast}(K;\mathfrak{g}_{\mathrm{Ad}_
\rho})$ is defined as the following alternating product
 $$\mathbb{T}(C_{\ast}(K;\mathfrak{g}_{\mathrm{Ad}_
\rho}),\{\mathbf{c}_p\}_{p=0}^{n},\{\mathbf{h}_p\}_{p=0}^{n})
 =\prod_{p=0}^n \left[\mathbf{b}_p\sqcup \ell_p(\mathbf{h}_p)\sqcup
s_p(\mathbf{b}_{p-1}), \mathbf{c}_p\right]^{(-1)^{(p+1)}},$$ where
$\left[\mathbf{e}_p, \mathbf{f}_p\right]$ is the determinant of
the transition matrix from basis $\mathbf{f}_p$ to $\mathbf{e}_p$
of $C_{p}.$
\end{definition}

By \cite{Porti2}, $\mathbb{T}(C_{\ast}(K;\mathfrak{g}_{\mathrm{Ad}_
\rho}),\{\mathbf{c}^p\}_{p=0}^{n},\{\mathbf{h}^p\}_{p=0}^{n})$ depends on the conjugacy class of $\rho$. It does not depend on a choice of the lifts of the cells in the universal cover of $X$ since the adjoint map is unimodular. It does not depend on  the basis $\mathcal{A}$ if the Euler characteristic $\chi(X)$ vanishes. Following \cite{Milnor}, it is also independent of the bases $\mathbf{b}^p$ and the sections $s_p, \ell_p.$

By \cite{Milnor}, the adjoint Reidemeister torsion is invariant under subdivisions, hence it defines an invariant of manifolds with dimension less than or equal to $3.$ More precisely,

\begin{definition}
Let $M$ be a smooth compact $n$-manifold with a triangulation $K,$ where $n\leq 3.$ For a given representation $\rho:\pi_1(M)\ra G$ and given bases 
$\{\mathbf{h}_p\}_{p=0}^{n}$ of homologies, the adjoint Reidemeister torsion of $M$ can be defined as
$$\mathbb{T}_{\rho}(M,\{\mathbf{h}_p\}_{p=0}^{n})=\mathbb{T}(C_{\ast}(K;\mathfrak{g}_{\mathrm{Ad}_
\rho}),\{\mathbf{c}_p\}_{p=0}^{n},\{\mathbf{h}_p\}_{p=0}^{n}).$$ 
\end{definition}


Mayer-Vietoris sequence is one of the useful tools to compute the adjoint Reidemeister torsion. More precisely, 

\begin{theorem}\label{prt1}
Assume that $X$ is a compact CW-complex with subcomplexes $X_1,$ $X_2 \subset X$ so that $X=X_1\cup X_2$ and $Y=X_1\cap X_2.$ Let $Y_1,\ldots,Y_k$ be the connected components of $Y.$ For $\nu= 1,2,$ consider the inclusions 
$$Y \overset{i_\nu} {\hookrightarrow} X_\nu \overset{j_\nu} {\hookrightarrow}X.$$  For $\nu= 1,2,$ and $\mu=1,\ldots,k,$ let $\rho:\pi_1(X)\rightarrow G$ be a representation with the restrictions 
$\rho_{_{|_{X_{\nu}}}}:\pi_1(X_\nu)\rightarrow G,$ 
$\rho_{_{|_{Y_{\mu}}}}:\pi_1(Y_\mu)\rightarrow G.$ Then there is a Mayer-Vietoris long exact sequence in homology with twisted coefficients 
\begin{eqnarray}\label{es12es}
& &\mathcal{H}_{\ast}:\cdots  \longrightarrow \underset {{\mu}}{\oplus} 
H_i(Y_{\mu};\mathfrak{g}_{\mathrm{Ad}_{\rho_{_{|_{Y_{\mu}}}}}})\overset{{i}^\ast_1 \oplus {i}^\ast_2 }{\longrightarrow }H_i(X_1;\mathfrak{g}_{\mathrm{Ad}_{\rho_{_{|_{X_{1}}}}}})\oplus H_i(X_2;\mathfrak{g}_{\mathrm{Ad}_{\rho_{_{|_{X_{2}}}}}}) \nonumber \\ 
&&\quad \quad\quad \quad\quad\quad\quad \quad\quad\; \tikz\draw[->,rounded corners,nodes={asymmetrical rectangle}](5.70,0.4)--(5.70,0)--(1,0)--(1,-0.4)node[yshift=4.0ex,xshift=12.0ex] {${j}^\ast_1 - {j}^\ast_2$ };
 \nonumber\\ 
&&\quad\quad\quad\quad\quad  H_i(X;\mathfrak{g}_{\mathrm{Ad}_\rho})\longrightarrow  \underset {\mu}{\oplus} H_{i-1}(Y_{\mu};\mathfrak{g}_{\mathrm{Ad}_{\rho_{_{|_{Y_{\mu}}}}}})\overset{{i}^\ast_1 \oplus {i}^\ast_2 }{\longrightarrow } \cdots
 \end{eqnarray}
Choose a basis for each of these homology groups such as $\mathbf{h}^X_i$ for 
$H_i(X;\mathfrak{g}_{\mathrm{Ad}_\rho}),$ $\mathbf{h}_i^{X_1}$ for $H_i(X_1;\mathfrak{g}_{\mathrm{Ad}_{\rho_{_{|_{X_{1}}}}}}),$ $\mathbf{h}^{X_2}_i$ for $H_i(X_2;\mathfrak{g}_{\mathrm{Ad}_{\rho_{_{|_{X_{2}}}}}}),$ and $\mathbf{h}^{Y_{\mu}}_i$ for 
$H_i(Y_{\mu};\mathfrak{g}_{\mathrm{Ad}_{\rho_{_{|_{Y_{\mu}}}}}}).$  
The long exact sequence (\ref{es12es}) can be viewed as a chain complex. Indeed, the following formula holds
\begin{eqnarray*}
&& \mathbb{T}_{\rho_{_{|_{X_{1}}}}}(X_1,\{\mathbf{h}^{X_1}_i\})\; 
 \mathbb{T}_{\rho_{_{|_{X_{2}}}}}(X_2,\{\mathbf{h}^{X_2}_i\})=\mathbb{T}_{\rho}(X,\{\mathbf{h}^X_i\})
\prod_{\mu=1}^k\mathbb{T}_{\rho_{_{|_{Y_{\mu}}}}}(Y_{\mu},\{\mathbf{h}^{Y_{\mu}}_i\})\\
&& \quad \quad \quad \quad \quad \quad  \quad\quad \quad \quad \quad  \quad   \quad \quad \quad \quad  \quad  \times \; \mathbb{T}(\mathcal{H}_{\ast}, \{\mathbf{h}_{\ast \ast}\}).
\end{eqnarray*}
\end{theorem} 
The proof of Theorem \ref{prt1} can be found in \cite{Milnor} or in \cite{Porti2}. For further information and the detailed proof, the reader is also referred to \cite{Milnor2,Porti,Turaev,Witten} and the references therein.
\begin{notation}
The torsion $\mathbb{T}(\mathcal{H}^{\ast},\{\mathbf{h}^{\ast \ast}\})$ is called the \textit{corrective term} and the same terminology is also used in \cite{BorgStefa}.
\end{notation}
\section{Main Results}
\subsection{Auxilary Results}
To prove Theorem \ref{theo1}, we give the following auxiliary results.

Let $\mathbb{D}^2$ be a closed disk. Since $\pi_1(\mathbb{D}^2)$ is trivial, any representation $\rho:\pi_1(\mathbb{D}^2)\rightarrow G$ is trivial. Note that the adjoint Reidemeister torsion is a simple homotopy invariant and $\mathbb{D}^2$ is simple homotopy equivalent to the point $*.$ Let $K=\tilde{K}=\{e_0\}$ denote the single $0$-cell of the point $*.$ As $\mathbb{Z}[\pi_1(\mathbb{D}^2)]=\mathbb{Z}[<e_0>]\cong \mathbb{Z},$
 \begin{equation}\label{defchn1}
 C_{0}(K;\mathrm{Ad}_{\rho})=C_{0}(\widetilde{K}; \mathbb{Z})\otimes 
 \mathfrak{g}/\sim \;\cong \mathfrak{g}.
 \end{equation}
 
Fix a $\mathcal{B}$-orthonormal $\mathbb{C}$-basis $\mathcal{A}=\{\mathfrak{a}_i\}_{i=1}^{\dim \mathfrak{g}}$ of $\mathfrak{g}.$ Then
$\mathbf{c}_0=\{e_0\otimes\mathfrak{a}_k\}_{k=1}^{\dim \mathfrak{g}}$
is the geometric basis of $C_{0}(K;{\mathrm{Ad}_\rho}).$ 

Consider the following chain complex
\begin{equation}\label{chaincomplex23}
  \begin{array}{ccc}
 C_{\ast}:= C_{\ast}(K;\mathrm{Ad}_{\rho})=(0 \stackrel{\partial_1}{\longrightarrow}  C_{0}(K;{\mathrm{Ad}_
\rho}) \stackrel{\partial_0}{\longrightarrow} 0).
  \end{array}
\end{equation}
Since the following equalities hold
\begin{eqnarray*}
&& B_0(C_{\ast})=\mathrm{Im}\{\partial_{1}:C_{1}( C_{\ast})\rightarrow C_{0}( C_{\ast})
\}=0,\\
&& Z_0(C_{\ast})=\mathrm{Ker}\{\partial_{0}:C_{0}( C_{\ast})\rightarrow
C_{-1}( C_{\ast})\}=C_{0}(K;{\mathrm{Ad}_\rho}),
\end{eqnarray*}
 the homology of $\mathbb{D}^2$ twisted by $\rho$ can be given as follows
$$H_0(\mathbb{D}^2;\mathrm{Ad}_{\rho})=H_0(C_{\ast})=Z_0(C_{\ast})/B_0(C_{\ast})\cong C_{0}(K;{\mathrm{Ad}_\rho})\cong \mathfrak{g}.$$ 

Then there are
the following short exact sequences
\begin{equation}\label{Equation1eks}
0\longrightarrow Z_0(C_\ast) \stackrel{\imath}{\hookrightarrow} C_0(C_{\ast})
\stackrel{\partial_0}{\longrightarrow} B_{-1}(C_\ast) \longrightarrow 0,
\end{equation}
\begin{equation}\label{Equation2eks}
0\longrightarrow B_0(C_\ast) \stackrel{\imath}{\hookrightarrow} Z_0(C_\ast)
\stackrel{\varphi_0}{\longrightarrow} H_0(C_\ast) \longrightarrow 0.
\end{equation}
Here, $\imath$ and $\varphi_p$ are the inclusion and the natural
projection, respectively.

Let $s_0:B_{-1}(C_{\ast})\rightarrow C_0(C_{\ast})$ and
$\ell_0:H_0(C_{\ast})\rightarrow Z_0(C_{\ast})$ be sections of the homomorphisms
$\partial_0:C_0(C_{\ast}) \rightarrow B_{-1}(C_{\ast}),$
$\varphi_0:Z_0(C_{\ast})\rightarrow H_0(C_{\ast}),$
respectively. Since $B_{0}(C_{\ast})=B_{-1}(C_{\ast})=\{0\},$ the homomorphism 
$\varphi_0$ becomes an isomorphism so the section $\ell_0$ is the inverse of this isomorphism. Moreover, we have
\begin{equation}\label{klur2eks}
C_0(C_{\ast})= \ell_0 (H_0(C_{\ast})).
\end{equation}
Let $\mathbf{h}^{\mathbb{D}^2}_0$ be an arbitrary basis of $H_0(\mathbb{D}^2;\mathrm{Ad}_{\rho}).$ From equation (\ref{klur2eks}) it follows
\begin{equation}\label{dısk1}
\mathbb{T}_{\rho}(\mathbb{D}^2,\{\mathbf{h}_0^{\mathbb{D}^2}\})=\left[\ell_0(\mathbf{h}_0^{\mathbb{D}^2}), \mathbf{c}_0\right].
\end{equation}

The following lemma is evident from equation (\ref{dısk1}).
\begin{lemma}\label{rem1dsk}
If we take the bases $\mathbf{h}_{0}^{\mathbb{D}^2}$ of 
$H_0(\mathbb{D}^2;\mathrm{Ad}_{\rho})$ as 
 $\varphi_0(\mathbf{c}_0)$ in equation~(\ref{dısk1}), then we get
\begin{equation*}
\mathbb{T}_{\rho}(\mathbb{D}^2,
\{\mathbf{h}^{\mathbb{D}^2}_0\})=\left[\ell_0(\mathbf{h}_0^{\mathbb{D}^2}), 
\mathbf{c}_0\right]=\left[\ell_0(\varphi_0(\mathbf{c}_0)),\mathbf{c}_0\right]=\left[\mathbf{c}_0,\mathbf{c}_0\right]=1.
\end{equation*}
\end{lemma}

\subsection{The Proof of Theorem \ref{theo1}}
Let $G$ be a complex reductive algebraic group $SL_n(\mathbb{C})$ or $PSL_n(\mathbb{C}).$ Let $M_1, M_2 \in \chi.$ We first focus on the case where $M$ is the disk sum of $M_1$ and $M_2$
$$M=M_1 \bigtriangleup M_2.$$ Cleary, $M \in \chi.$ From the Seifert-Van Kampen's theorem it follows
$$\pi_1(M) = \pi_1(M_1)\ast \pi_1(M_2).$$

Let $\psi_{_1}:\pi_1(M_1)\rightarrow G$ and $\psi_{_2}:\pi_1(M_1)\rightarrow G$  be representations. By the universal property of free product, there is a unique homomorphism 
$$\varrho:\pi_1(M_1)\ast \pi_1(M_2)\rightarrow G$$ with the restrictions 
$\varrho_{|_{\pi_1(M_1)}}=\psi_{_1}$ and $\varrho_{|_{\pi_1(M_2)}}=\psi_{_2}.$
Moreover, we consider the restriction $\varrho_{|_{\mathbb{D}^2}}$ of the representation $\varrho$ to $\pi_1({\mathbb{D}^2}).$ Then there is a short exact sequence of the chain complexes
\begin{small}
\begin{equation}\label{seq2PP}
0\to C_{\ast}(\mathbb{D}^2;\mathfrak{g}_{\mathrm{Ad}_{\varrho_{|_{\mathbb{D}^2}}}})\rightarrow 
C_{\ast}(M_1;\mathfrak{g}_{\mathrm{Ad}_{\psi_{_1}}})\oplus
C_{\ast}(M_2;\mathfrak{g}_{\mathrm{Ad}_{\psi_{_2}}}) \rightarrow C_{\ast}(M;\mathfrak{g}_{\mathrm{Ad}_\varrho})\to 0.
\end{equation}
\end{small}

It is well–known that for any $3$–manifold $N$ with non-empty boundary, one has $\mathrm{dim}(H_1(N;\mathfrak{g}_{\mathrm{Ad}_\rho}))\geq 1,$ where 
$\rho:\pi_1(N)\rightarrow G$ is a representation. Thus, $C_{\ast}(N;\mathfrak{g}_{\mathrm{Ad}_\rho})$ is never acyclic. Hence, associated to the short exact sequence (\ref{seq2PP}), by \cite{Porti2}, there is a Mayer-Vietoris long exact sequence in homology with twisted coefficients
\begin{small}
\begin{eqnarray*}\label{longexact1PP}
&\mathcal{H}_{\ast}:&0 \stackrel{\partial'_3}{\rightarrow} H_3(M_1;\mathfrak{g}_{\mathrm{Ad}_{\psi_{_1}}})\oplus H_3(M_2;\mathfrak{g}_{\mathrm{Ad}_{\psi_{_2}}})
\stackrel{\partial_3}{\rightarrow} H_3(M;\mathfrak{g}_{\mathrm{Ad}_\varrho})\\ \nonumber
&&\quad \quad \quad \tikz\draw[->,rounded corners,nodes={asymmetrical rectangle}](6.40,0.4)--(6.40,0)--(1,0)--(1,-0.4) node[yshift=4.5ex,xshift=16.0ex] {$\varphi_3$}; \\ \nonumber
& &0 \stackrel{\partial'_2}{\rightarrow} H_2(M_1;\mathfrak{g}_{\mathrm{Ad}_{\psi_{_1}}})\oplus H_2(M_2;\mathfrak{g}_{\mathrm{Ad}_{\psi_{_2}}})
\stackrel{\partial_2}{\rightarrow} H_2(M;\mathfrak{g}_{\mathrm{Ad}_\varrho}) \\ \nonumber
 &&\quad \quad \quad \tikz\draw[->,rounded corners](6.40,0.4)--(6.40,0)--(1,0)--(1,-0.4) node[yshift=4.5ex,xshift=16.0ex] {$\varphi_2$}; \\  \nonumber
& &0 \stackrel{\partial'_1}{\rightarrow} H_1(M_1;\mathfrak{g}_{\mathrm{Ad}_{\psi_{_1}}})\oplus H_1(M_2;\mathfrak{g}_{\mathrm{Ad}_{\psi_{_2}}})
\stackrel{\partial_1}{\rightarrow} H_1(M;\mathfrak{g}_{\mathrm{Ad}_\varrho}) \\
&&\quad \quad \quad \tikz\draw[->,rounded corners](6.40,0.4)--(6.40,0)--(1,0)--(1,-0.4) node[yshift=4.5ex,xshift=16.0ex] {$\varphi_1$}; \\  \nonumber
& &H_0(\mathbb{D}^2;\mathfrak{g}_{\mathrm{Ad}_{\varrho_{|_{\mathbb{D}^2}}}})\stackrel{\partial'_0}{\rightarrow} H_0(M_1;\mathfrak{g}_{\mathrm{Ad}_{\psi_{_1}}})\oplus H_0(M_2;\mathfrak{g}_{\mathrm{Ad}_{\psi_{_2}}})
\stackrel{\partial_0}{\rightarrow} H_0(M;\mathfrak{g}_{\mathrm{Ad}_\varrho}) \stackrel{\varphi_0}{\rightarrow} 0.
  \end{eqnarray*}
\end{small}

\noindent From the exactness of $\mathcal{H}_{\ast}$ and the First Isomorphism Theorem it follows that $\varphi_1$ is a zero map and for each 
$i\in \{1,2,3\},$ the following isomorphism holds
$$ H_i(M_1;\mathfrak{g}_{\mathrm{Ad}_{\psi_{_1}}})\oplus H_i(M_2;\mathfrak{g}_{\mathrm{Ad}_{\psi_{_2}}})\overset {\partial_{i}}{\cong}  H_i(M;\mathfrak{g}_{\mathrm{Ad}_\varrho}).$$  

 Let $\mathbf{h}^{M}_p$ and $\mathbf{h}^{\mathbb{D}^2}_0$ be given bases of $H_p(M;\mathfrak{g}_{\mathrm{Ad}_\varrho})$ and $H_0(\mathbb{D}^2;\mathfrak{g}_{\mathrm{Ad}_{\varrho_{|_{\mathbb{D}^2}}}}).$ First, we denote the vector spaces in the long exact sequence $\mathcal{H}_{\ast}$ (from right to left) by $C_p(\mathcal{H}_{\ast})$ for $p \in \{0,\ldots,11\}.$ Then we consider the short exact sequences 
\begin{equation}\label{x1s}
0\to Z_p(\mathcal{H}_{\ast}) \hookrightarrow C_p(\mathcal{H}_{\ast})\overset{\partial_p}{\rightarrow} B_{p-1}(\mathcal{H}_{\ast}) \to 0,
\end{equation}
\begin{equation}\label{x2s} 
0\to B_p(\mathcal{H}_{\ast}) \hookrightarrow Z_p(\mathcal{H}_{\ast})\overset{\varphi_p}{\twoheadrightarrow} H_p(\mathcal{H}_{\ast}) \to 0.
\end{equation}
 For each $p,$ let us consider the isomorphism $s_{p}:B_{p-1}(\mathcal{H}_{\ast}) \rightarrow s_{p}(B_{p-1}(\mathcal{H}_{\ast}))$ obtained by the First Isomorphism Theorem as a section of $C_p(\mathcal{H}_{\ast})\rightarrow B_{p-1}(\mathcal{H}_{\ast}).$ Using the exactness of $\mathcal{H}_{\ast}$ in the short exact sequence (\ref{x2s}), we obtain
 $$B_p(\mathcal{H}_{\ast})=Z_p(\mathcal{H}_{\ast}).$$
  Hence, the sequence (\ref{x1s}) becomes
\begin{equation}\label{kfrdl1}
0\to B_p(\mathcal{H}_{\ast}) \hookrightarrow C_p(\mathcal{H}_{\ast})
 \twoheadrightarrow B_{p-1}(\mathcal{H}_{\ast}) \to 0.
\end{equation}
Applying the Splitting Lemma for the sequence (\ref{kfrdl1}), we have
\begin{equation}\label{uzuntamdzm7}
C_p(\mathcal{H}_{\ast})=B_p(\mathcal{H}_{\ast})\oplus s_{_{p}}(B_{p-1}(\mathcal{H}_{\ast})).
\end{equation}
Thus, the Reidemeister torsion of $\mathcal{H}_{\ast}$ satisfies the following formula
\begin{equation}\label{torlg}
\mathbb{T}\left( \mathcal{H}_{\ast},\{\mathbf{h}_p\}_{p=0}^{11},\{0\}_{p=0}^{11}\right)
 =\prod_{p=0}^{11} \left[\mathbf{h}'_p, \mathbf{h}_p\right]^{(-1)^{(p+1)}}.
\end{equation}
Here, $\mathbf{h}'_p$ denotes the new basis $\mathbf{b}_p\sqcup s_p(\mathbf{b}_{p-1})$ of $C_p(\mathcal{H}_{\ast})$ for all $p.$ Since Reidemeister torsion does not depend on bases $\mathbf{b}_p$ and sections $s_p,$ we can chose the suitable bases $\mathbf{b}_p$ and sections $s_p$ to show that the existence of the homology bases $\mathbf{h}^{M_1}_p,$ $\mathbf{h}^{M_2}_p$ in which the Reidemeister torsion of $\mathcal{H}_{\ast}$ in equation (\ref{torlg}) is equal to $1.$ 

$\bullet$ We consider the first part of the long exact sequence 
$\mathcal{H}_{\ast}:$
\begin{small}
 $$0\stackrel{\varphi_1}{\rightarrow}H_0(\mathbb{D}^2;\mathfrak{g}_{\mathrm{Ad}_{\varrho_{|_{\mathbb{D}^2}}}})\stackrel{\partial'_0}{\rightarrow} H_0(M_1;\mathfrak{g}_{\mathrm{Ad}_{\psi_{_1}}})\oplus H_0(M_2;\mathfrak{g}_{\mathrm{Ad}_{\psi_{_2}}})
\stackrel{\partial_0}{\rightarrow} H_0(M;\mathfrak{g}_{\mathrm{Ad}_\varrho}) \stackrel{\varphi_0}{\rightarrow} 0.$$
\end{small}

 First, we use equation (\ref{uzuntamdzm7}) for the vector space $C_0(\mathcal{H}_{\ast})=H_0(M;\mathfrak{g}_{\mathrm{Ad}_\varrho}).$ Since $\mathrm{Im}\,\varphi_0=\{0\},$ we get
\begin{equation}\label{estlk5677}
C_0(\mathcal{H}_{\ast})=\mathrm{Im}\,\partial_0\oplus s_{_0}(\mathrm{Im}\,\varphi_0)=\mathrm{Im}\,\partial_0.
\end{equation}
If we choose the basis $\mathbf{h}^{\mathrm{Im}\,\partial_0}$ of $\mathrm{Im}\,\partial_0$ as $\mathbf{h}_0^{M},$ then $\mathbf{h}_0^{M}$ becomes the new basis $\mathbf{h}'_0$ of $C_0(\mathcal{H}_{\ast})$ by equation (\ref{estlk5677}). Since
$\mathbf{h}_0^{M}$ is also the given basis $\mathbf{h}_0$ of $C_0(\mathcal{H}_{\ast}),$ the following equality holds
\begin{equation}\label{cnmbr07}
[\mathbf{h}'_0,\mathbf{h}_0]=1. 
\end{equation}

Let us consider equation (\ref{uzuntamdzm7}) for 
$C_1(\mathcal{H}_{\ast})=H_0(M_1;\mathfrak{g}_{\mathrm{Ad}_{\psi_{_1}}})\oplus H_0(M_2;\mathfrak{g}_{\mathrm{Ad}_{\psi_{_2}}}).$ Then we obtain
\begin{equation}\label{estlk1234s7} 
C_1(\mathcal{H}_{\ast})=\mathrm{Im}\,\partial'_0\oplus s_{_{1}}(\mathrm{Im}\,\partial_0).
\end{equation}
Recall that the basis $\mathbf{h}^{\mathrm{Im}\,\partial_0}$ of $\mathrm{Im}\,\partial_0$ was chosen as $\mathbf{h}_0^{M}$ in the previous step. Note also that $\mathrm{Im}\,\partial'_0$ is isomorphic to $H_0(\mathbb{D}^2;\mathfrak{g}_{\mathrm{Ad}_{\varrho_{|_{\mathbb{D}^2}}}})$ so we can take the basis 
$\mathbf{h}^{\mathrm{Im}\,\partial'_0}$ of $\mathrm{Im}\,\partial'_0$ as $\partial'_0(\mathbf{h}_0^{\mathbb{D}^2}).$ By equation (\ref{estlk1234s7}), we get a new basis for $C_1(\mathcal{H}_{\ast})$ as follows
$$\mathbf{h}'_1=\left\{\partial'_0(\mathbf{h}_0^{\mathbb{D}^2}),
s_{_{1}}(\mathbf{h}_0^{M})\right\}.$$ 

Let $n_0^{M},$ $n_0^{M_{\ell}},$ and $n_0^{\mathbb{D}^2}$ denote the dimension of spaces $H_0(M;\mathfrak{g}_{\mathrm{Ad}_\varrho}),$ $H_0(M_{\ell};\mathfrak{g}_{\mathrm{Ad}_{\psi_{_{\ell}}}}),$ and $H_0(\mathbb{D}^2;\mathfrak{g}_{\mathrm{Ad}_{\varrho_{|_{\mathbb{D}^2}}}})$ for $\ell=1,2.$ Since $H_0(M_1;\mathfrak{g}_{\mathrm{Ad}_{\psi_{_1}}})$ and $H_0(M_2;\mathfrak{g}_{\mathrm{Ad}_{\psi_{_2}}})$ are subspaces of $C_1(\mathcal{H}_{\ast}),$ we have
 \begin{eqnarray*}
n_0^{M_1}+n_0^{M_2}=n_0^M+n_0^{\mathbb{D}^2}=\mathrm{dim}(C_{1}(\mathcal{H}_{\ast})).
 \end{eqnarray*}
Recall that $\mathbf{h}_0^{M}=\left\{\mathbf{h}_{_{0,j}}^{M}\right\}_{j=1}^{n_0^M}$ is the given basis of $H_0(M;\mathfrak{g}_{\mathrm{Ad}_\varrho}).$ There are non-zero vectors
$(a_{_{i,1}},a_{_{i,2}},\cdots,a_{_{i,n_0^{M_1}+n_0^{M_2}}})$ for $i \in \{1,2,\ldots,n_0^{M_1}+n_0^{M_2}\}$ such that
\begin{eqnarray*}
& &\left\{\sum_{j=1}^{n_0^{\mathbb{D}^2}}a_{_{i,j}}\partial'_0(\mathbf{h}_0^{\mathbb{D}^2})+
\sum_{j=n_0^{\mathbb{D}^2}+1}^{n_0^M+n_0^{\mathbb{D}^2}}a_{_{i,j}}s_{_{1}}\left(\mathbf{h}_{_{0,j}}^{M}\right)\right\}_{i=1}^{n_0^{M_1}},\\
& &\left\{\sum_{j=1}^{n_0^{\mathbb{D}^2}}a_{_{i,j}}\partial'_0(\mathbf{h}_0^{\mathbb{D}^2})+
\sum_{j=n_0^{\mathbb{D}^2}+1}^{n_0^M+n_0^{\mathbb{D}^2}}a_{_{i,j}}s_{_{1}}\left(\mathbf{h}_{_{0,j}}^{M}\right)\right\}_{i=n_0^{M_1}+1}^{n_0^{M_1}+n_0^{M_2}}
\end{eqnarray*}
are bases of $H_0(M_1;\mathfrak{g}_{\mathrm{Ad}_{\psi_{_1}}})$ and $H_0(M_2;\mathfrak{g}_{\mathrm{Ad}_{\psi_{_2}}}),$ respectively. Moreover, 
$A=[a_{_{i,j}}]$ is the $((n_0^{M_1}+n_0^{M_2})\times (n_0^{M_1}+n_0^{M_2}))$ invertible matrix.

 Let us take the bases of $H_0(M_1;\mathfrak{g}_{\mathrm{Ad}_{\psi_{_1}}})$ and $H_0(M_2;\mathfrak{g}_{\mathrm{Ad}_{\psi_{_2}}})$ as follows
\begin{eqnarray*}
&&\mathbf{h}_0^{M_1}=\left\{(\det A)^{-1}\left[\sum_{j=1}^{n_0^{\mathbb{D}^2}}a_{_{1,j}}\partial'_0(\mathbf{h}_0^{\mathbb{D}^2})+\sum_{j=n_0^{\mathbb{D}^2}+1}^{n_0^M+n_0^{\mathbb{D}^2}}a_{_{1,j}}s_{_{1}}\left(\mathbf{h}_{_{0,j}}^{M}\right)\right], \right.\\
& & \quad \quad \quad \quad\quad\quad \quad\quad\; \left. \left\{\sum_{j=1}^{n_0^{\mathbb{D}^2}}a_{_{i,j}}\partial'_0(\mathbf{h}_0^{\mathbb{D}^2})+
\sum_{j=n_0^{\mathbb{D}^2}+1}^{n_0^M+n_0^{\mathbb{D}^2}}a_{_{i,j}}s_{_{1}}\left(\mathbf{h}_{_{0,j}}^{M}\right)\right\}_{i=2}^{n_0^{M_1}}\right\},\\
&& \mathbf{h}_0^{M_2}=
 \left\{\sum_{j=1}^{n_0^{\mathbb{D}^2}}a_{_{i,j}}\partial'_0(\mathbf{h}_0^{\mathbb{D}^2})+
\sum_{j=n_0^{\mathbb{D}^2}+1}^{n_0^M+n_0^{\mathbb{D}^2}}a_{_{i,j}}s_{_{1}}\left(\mathbf{h}_{_{0,j}}^{M}\right)\right\}_{i=n_0^{M_1}+1}^{n_0^{M_1}+n_0^{M_2}}.
\end{eqnarray*}
\noindent Then $\mathbf{h}_1=\{\mathbf{h}_0^{M_1}, \mathbf{h}_0^{M_2}\}$ 
 becomes the initial basis of $C_1(\mathcal{H}_{\ast})$ and we have
\begin{equation}\label{cnmbr17}
[\mathbf{h}'_1,\mathbf{h}_1]=1.
\end{equation}

 Considering the space $C_2(\mathcal{H}_{\ast})=H_0(\mathbb{D}^2;\mathfrak{g}_{\mathrm{Ad}_{\varrho_{|_{\mathbb{D}^2}}}})$ in equation (\ref{uzuntamdzm7}) and using the fact that 
$\mathrm{Im}\,\varphi_1=\{0\},$ we can write the space $C_2(\mathcal{H}_{\ast})$ as follows
\begin{equation}\label{estlk3k7} 
C_2(\mathcal{H}_{\ast})=\mathrm{Im}\,\varphi_1\oplus s_{_{2}}(\mathrm{Im}\,\partial'_0)= s_{_{2}}(\mathrm{Im}\,\partial'_0).
\end{equation}
 By equation (\ref{estlk3k7}), $s_{_{2}}(\partial'_0(\mathbf{h}_0^{\mathbb{D}^2}))=\mathbf{h}_0^{\mathbb{D}^2}$ becomes a new basis 
$\mathbf{h}'_2$ of $C_2(\mathcal{H}_{\ast}).$ Since 
$\mathbf{h}_0^{\mathbb{D}^2}$ is also the given basis $\mathbf{h}_2$ of $C_2(\mathcal{H}_{\ast}),$ we get
\begin{equation}\label{cnmbr27}
[\mathbf{h}'_2,\mathbf{h}_2] =1. 
\end{equation}

For the space $C_3(\mathcal{H}_{\ast})=\{0\}$ in the sequence 
$\mathcal{H}_{\ast},$ we use the convention $1\cdot0=1.$ Then we obtain
\begin{equation}\label{cnmbr37}
[\mathbf{h}'_3,\mathbf{h}_3] =1.
\end{equation}

$\bullet$  Let us consider the second part of the sequence $\mathcal{H}_{\ast}$ for $i=1,2,3$ 
 \begin{equation}\label{eq1}
 0 \stackrel{\partial'_i}{\rightarrow} H_i(M_1;\mathfrak{g}_{\mathrm{Ad}_{\psi_{_1}}})\oplus H_i(M_2;\mathfrak{g}_{\mathrm{Ad}_{\psi_{_2}}})
\stackrel{\partial_i}{\rightarrow} H_i(M;\mathfrak{g}_{\mathrm{Ad}_\varrho})\stackrel{\partial''_i}{\rightarrow} 0.
 \end{equation}

 Now we denote the vector spaces in the short exact sequence (\ref{eq1}) (from right to left) as $C_{3i}(\mathcal{H}_{\ast}),$ $C_{3i+1}(\mathcal{H}_{\ast})$ and $C_{3i+2}(\mathcal{H}_{\ast})$ for each $i=1,2,3.$

Note that the spaces $C_{3i}(\mathcal{H}_{\ast})$ equals to $\{0\}.$ If we use the convention $1\cdot0=1$ for each $i=2,3.$ Then we get
\begin{equation}\label{cnvtyr37}
[\mathbf{h}'_{3i},\mathbf{h}_{3i}] =1.
\end{equation}

 By the exactness of $\mathcal{H}_{\ast},$ we get the following isomorphism:
$$H_i(M_1;\mathfrak{g}_{\mathrm{Ad}_{\psi_{_1}}})\oplus H_i(M_2;\mathfrak{g}_{\mathrm{Ad}_{\psi_{_2}}})\overset {\partial_{i}}{\cong}  H_i(M;\mathfrak{g}_{\mathrm{Ad}_\varrho})$$  
for each $i=1,2,3.$

 We use equation (\ref{uzuntamdzm7}) for the space
$C_{{3i+1}}(\mathcal{H}_{\ast})=H_i(M;\mathfrak{g}_{\mathrm{Ad}_\varrho}).$ 
Since $\mathrm{Im}\,\partial''_i=\{0\},$ the following equality holds
\begin{equation}\label{estlk12m7} 
C_{{3i+1}}(\mathcal{H}_{\ast})=\mathrm{Im}\,\partial_i\oplus s_{_{3i+1}}(\mathrm{Im}\,\partial''_i)=\mathrm{Im}\,\partial_i.
\end{equation}
Since ${\mathrm{Im}\,\partial_i}$ equals to $H_i(M;\mathfrak{g}_{\mathrm{Ad}_\varrho}),$ we can take the basis $\mathbf{h}^{\mathrm{Im}\,\partial_i}$ of ${\mathrm{Im}\,\partial_i}$ as $\mathbf{h}_i^{M}.$ By equation (\ref{estlk12m7}), $\mathbf{h}_i^{M}$ becomes a new basis $\mathbf{h}'_{{3i+1}}$ of 
$C_{{3i+1}}(\mathcal{H}_{\ast}).$ As $\mathbf{h}_i^{M}$ is also the given basis $\mathbf{h}_{{3i+1}}$ of $C_{{3i+1}}(\mathcal{H}_{\ast}),$ the following equation holds
\begin{equation}\label{cnmbr47}
 [\mathbf{h}'_{{3i+1}},\mathbf{h}_{{3i+1}}] =1. 
 \end{equation}
 
 Considering equation (\ref{uzuntamdzm7}) for
$C_{{3i+2}}(\mathcal{H}_{\ast})=H_i(M_1;\mathfrak{g}_{\mathrm{Ad}_{\psi_{_1}}})\oplus H_i(M_2;\mathfrak{g}_{\mathrm{Ad}_{\psi_{_2}}})$ and using the fact that $\mathrm{Im}\,\partial'_i=\{0\},$ we obtain
\begin{equation}\label{estytrn5} 
C_{{3i+2}}(\mathcal{H}_{\ast})=\mathrm{Im}\,\partial'_i\oplus s_{_{3i+2}}(\mathrm{Im}\,\partial_i)=s_{_{3i+2}}(\mathrm{Im}\,\partial_i).
\end{equation}
Since $H_i(M_1;\mathfrak{g}_{\mathrm{Ad}_{\psi_{_1}}})\oplus H_i(M_2;\mathfrak{g}_{\mathrm{Ad}_{\psi_{_2}}})$ and $H_i(M;\mathfrak{g}_{\mathrm{Ad}_\varrho})$ are isomorphic, the section $s_{_{3i+2}}$ can be considered as the inverse of the isomorphism $\partial_i.$ In the previous step, the basis $\mathbf{h}^{\mathrm{Im}\,\partial_i}$ of ${\mathrm{Im}\,\partial_i}$ was chosen as $\mathbf{h}_i^{M}.$ By equation (\ref{estytrn5}), 
$s_{_{3i+2}}(\mathbf{h}_i^{M})$ becomes a new basis
$\mathbf{h}'_{3i+2}$ of $C_{{3i+2}}(\mathcal{H}_{\ast}).$

 Let $n_i^{M}$ and $n_i^{M_{\ell}}$ denote the dimension of spaces $H_i(M;\mathfrak{g}_{\mathrm{Ad}_\varrho})$ and $H_i(M_{\ell};\mathfrak{g}_{\mathrm{Ad}_{\psi_{_{\ell}}}})$ for $\ell=1,2.$  Note that $H_i(M_1;\mathfrak{g}_{\mathrm{Ad}_{\psi_{_1}}})$ and $H_i(M_2;\mathfrak{g}_{\mathrm{Ad}_{\psi_{_2}}})$ are subspaces of $C_{3i+2}(\mathcal{H}_{\ast})$ and $\mathbf{h}_i^{M}=\left\{\mathbf{h}_{_{i,j}}^{M}\right\}_{j=1}^{n_i^{M}}$ is the given basis of $H_i(M;\mathfrak{g}_{\mathrm{Ad}_\varrho}).$ Moreover, the exactness of $\mathcal{H}_{\ast}$ yields
 \begin{eqnarray*}
n_i^{M_1}+n_i^{M_2}=n_i^M=\mathrm{dim}(C_{3i+2}(\mathcal{H}_{\ast})).
 \end{eqnarray*}
There are non-zero vectors
$(a_{_{k,1}},a_{_{k,2}},\cdots,a_{_k,n_i^M})$ for $k\in \{1,2,\ldots,n_i^{M}\}$ such that
\begin{eqnarray*}
\left\{
\sum_{j=1}^{n_i^{M}}a_{_{k,j}}s_{_{3i+2}}\left(\mathbf{h}_{_{i,j}}^{M}\right)\right\}_{k=1}^{n_i^{M_1}},\left\{
\sum_{j=1}^{n_i^{M}}a_{_{k,j}}s_{_{3i+2}}\left(\mathbf{h}_{_{i,j}}^{M}\right)\right\}_{k=n_i^{M_1}+1}^{n_i^{M}}
\end{eqnarray*}
are bases of $H_i(M_1;\mathfrak{g}_{\mathrm{Ad}_{\psi_{_1}}})$ and $H_i(M_2;\mathfrak{g}_{\mathrm{Ad}_{\psi_{_2}}}),$ respectively. Moreover, 
$A=[a_{_{i,j}}]$ is the $(n_i^M\times n_i^M)$ invertible matrix .

 Let us take the bases of $H_i(M_1;\mathfrak{g}_{\mathrm{Ad}_{\psi_{_1}}})$ and $H_i(M_2;\mathfrak{g}_{\mathrm{Ad}_{\psi_{_2}}})$ as follows
\begin{eqnarray*}
&&\mathbf{h}_{3i+2}^{M_1}=\left\{(\det A)^{-1}\left[\sum_{j=1}^{n_i^M}a_{_{1,j}}s_{_{3i+2}}\left(\mathbf{h}_{_{i,j}}^{M}\right)\right], \left\{
\sum_{j=1}^{n_i^M}a_{_{k,j}}s_{_{3i+2}}\left(\mathbf{h}_{_{i,j}}^{M}\right)\right\}_{k=2}^{n_i^{M_1}}\right\}\\
&& \mathbf{h}_{3i+2}^{M_2}=
 \left\{
\sum_{j=1}^{n_i^{M}}a_{_{k,j}}s_{_{3i+2}}\left(\mathbf{h}_{_{i,j}}^{M}\right)\right\}_{k=n_i^{M_1}+1}^{n_i^{M}}.
\end{eqnarray*}
Then $\mathbf{h}_{3i+2}=\{\mathbf{h}_{3i+2}^{M_1}, \mathbf{h}_{3i+2}^{M_2}\}$ becomes the initial basis of $C_{3i+2}(\mathcal{H}_{\ast})$ and we have
\begin{equation}\label{cnmbrxd4}
 [\mathbf{h}'_{{3i+2}},\mathbf{h}_{{3i+2}}] =1. 
 \end{equation}

 Equations (\ref{cnmbr07}), (\ref{cnmbr17}), 
(\ref{cnmbr27}), (\ref{cnmbr37}), (\ref{cnmbr47}), (\ref{cnmbrxd4}) yield
\begin{equation}\label{uznhm11tm7}
  \mathbb{T}( \mathcal{H}_{\ast},\{\mathbf{h}_p\}_{p=0}^{11} ,\{0\}_{p=0}^{11} )
  =\prod_{p=0}^{11} \left[\mathbf{h}'_p,\mathbf{h}_p\right]^{(-1)^{(p+1)}} =1. 
\end{equation}

Combining Theorem~\ref{prt1} and equation (\ref{uznhm11tm7}), the following formula is valid
\begin{equation}\label{eq:123}
\mathbb{T}_{\varrho}(M,\{\mathbf{h}_p^{M}\}_{p=0}^3)=\frac{
\mathbb{T}_{{\psi_{_1}}}(M_1,\{\mathbf{h}_p^{M_1}\}_{p=0}^{3}) 
\;\mathbb{T}_{{\psi_{_2}}}(M_2,\{\mathbf{h}_p^{M_2}\}_{p=0}^{3})}
{\mathbb{T}_{{\varrho_{|_{\mathbb{D}^2}}}}(\mathbb{D}^2,\{\mathbf{h}_0^{\mathbb{D}^2}\})}.
\end{equation}

By Lemma \ref{rem1dsk} and equation (\ref{eq:123}), we have
\begin{equation}\label{eq:12321}
\mathbb{T}_{\varrho}(M,\{\mathbf{h}_p^{M}\}_{p=0}^3)=\mathbb{T}_{{\psi_{_1}}}(M_1,\{\mathbf{h}_p^{M_1}\}_{p=0}^{3}) \;\mathbb{T}_{{\psi_{_2}}}(M_2,\{\mathbf{h}_p^{M_2}\}_{p=0}^{3}).
\end{equation}

Applying equation (\ref{eq:12321}) inductively finishes the proof of Theorem \ref{theo1}.


\begin{thebibliography}{33}

\bibitem{BorgStefa}
Borghini, S., A gluing formula for reidemeister--turaev torsion,
 Annali di Matematica Pura ed Applicata, \textbf{194} (5), 1535--1561, (2015)



\bibitem{T. A. Chap} Chapman, T. A.:Topological invariance of Whitehead torsion, Amer. J. Math. \textbf{96} 3, 488--497 (1974)


\bibitem{Franz} Franz W. : \"{U}ber die Torsion einer \"{U}berdeckung. J. Reine Angew. Math. \textbf{173}, 245--254 (1935) 

\bibitem{RCLC} Kirby R. C., Siebenmann L. C.: On triangulation of manifolds and Haupvermutung. Bull. Amer. Math. Soc. \textbf{75}, 742--749 (1969)

\bibitem{JLGROSS1} Gross J.L., A unique decomposition theorem for 3-manifolds with connected boundary, Trans. Amer. Math. Soc. 142 (1969), 191--199.



\bibitem{Milnor2} Milnor J.: A duality theorem for Reidemeister torsion. Ann. of Math. \textbf{76} (1), 137--147 (1962) 

\bibitem{Milnor} Milnor J.: Whitehead torsion. Bull. Amer. Math. Soc. 
\textbf{72}, 358--426 (1966)

\bibitem{Porti} Porti J.: Torsion de Reidemeister pour les Varieties Hyperboliques.
Memoirs of the Amer. Math. Soc. 128 (612), Amer. Math. Soc., Providence, RI (1997).

\bibitem{Porti2} Porti J.: Reidemeister torsion, hyperbolic three-manifolds, and character varieties, arXiv:1511.00400v2, (2016).



\bibitem{Reidemeister} Reidemeister K.: Homotopieringe und
Linsenr\"{a}ume. Hamburger Abhandl. \textbf{11}, 102--109 (1935)

\bibitem{Witten} Witten E.: On quantum gauge theories in twodimensions. Comm. Math. Phys. \textbf{141}, 153--209 (1991)


\bibitem{Turaev} Turaev V.: Torsions of 3-Dimensional Manifolds. \textbf{208}. Birkhauser Verlag (2002)

\end{thebibliography}
\end{document}